\documentclass[10.5pt]{article}
\usepackage{amsfonts,mathrsfs}

\def\sqr#1#2{{\vcenter{\vbox{\hrule height.#2pt
              \hbox{\vrule width.#2pt height#1pt \kern#1pt \vrule width.#2pt}
              \hrule height.#2pt}}}}
\def\signed #1{{\unskip\nobreak\hfil\penalty50
              \hskip2em\hbox{}\nobreak\hfil#1
              \parfillskip=0pt \finalhyphendemerits=0 \par}}
\def\endpf{\signed {$\sqr69$}}
\def\3n{\negthinspace \negthinspace \negthinspace }
\def\2n{\negthinspace \negthinspace }
\def\1n{\negthinspace }

\def\dbE{\mathbb{E}}
\def\dbF{\mathbb{F}}

\def\dbN{\mathbb{N}}

\def\dbR{\mathbb{R}}

\def\dbY{\mathbb{Y}}
\def\dbZ{\mathbb{Z}}

\def\sL{\mathscr{L}}

\def\sP{\mathscr{P}}

\def\={\buildrel \triangle \over =}

\def\ds{\displaystyle}

\def\ns{\noalign{\ss}}

\def\b{\beta}
\def\g{\gamma}
\def\d{\delta}
\def\e{\varepsilon}

\def\l{\lambda}

\def\si{\sigma}
\def\t{\tau}
\def\f{\varphi}

\def\Th{\Theta}
\def\L{\Lambda}

\def\O{\Omega}

\def\cF{{\cal F}}

\def\cL{{\cal L}}

\def\cP{{\cal P}}

\def\ss{\smallskip}
\def\ms{\medskip}

\def\q{\quad}
\def\qq{\qquad}
\def\hb{\hbox}

\def\esssup{\mathop{\rm esssup}}

\def\h{\widehat}

\def\cd{\cdot}

\def\({\Big (}
\def\){\Big )}
\def\[{\Big[}
\def\]{\Big]}
\def\bde{\begin{definition}}
\def\ede{\end{definition}}
\def\be{\begin{equation}}
\def\bel{\begin{equation}\label}
\def\ee{\end{equation}}
\def\bt{\begin{theorem}}
\def\et{\end{theorem}}
\def\bc{\begin{corollary}}
\def\ec{\end{corollary}}
\def\bl{\begin{lemma}}
\def\el{\end{lemma}}
\def\bp{\begin{proposition}}
\def\ep{\end{proposition}}
\def\bas{\begin{assumption}}
\def\eas{\end{assumption}}
\def\br{\begin{remark}}
\def\er{\end{remark}}
\def\ba{\begin{array}}
\def\ea{\end{array}}
\def\ed{\end{document}}

\def\square#1{\vbox{\hrule\hbox{\vrule height#1%
     \kern#1\vrule}\hrule}}
\def\rectangle#1#2{\vbox{\hrule\hbox{\vrule height#1%
     \kern#2\vrule}\hrule}}

\font\tenbb=msbm10 \font\sevenbb=msbm7 \font\fivebb=msbm5

\newfam\bbfam
\scriptscriptfont\bbfam=\fivebb \textfont\bbfam=\tenbb
\scriptfont\bbfam=\sevenbb

\newtheorem{lemma}{Lemma}[section]
\newtheorem{remark}{Remark}[section]

\newtheorem{theorem}{Theorem}[section]
\newtheorem{corollary}{Corollary}[section]

\newtheorem{definition}{Definition}[section]
\newtheorem{proposition}{Proposition}[section]
\newtheorem{assumption}{Assumption}[section]

\makeatletter
   
   \@addtoreset{equation}{section}
\makeatother

\usepackage{amsfonts}
\usepackage{latexsym}
\usepackage{amsmath}
\usepackage{amssymb}
\usepackage{color}

\begin{document}

\title{Mean-variance portfolio selection and variance hedging with random coefficients: closed-loop equilibrium strategy\footnote{The research was supported by the NSF of China under grant 11231007, 11401404 and 11471231.}}

\author{Tianxiao Wang \footnote{School of
Mathematics, Sichuan University, Chengdu, P. R. China. Email:wtxiao2014@scu.edu.cn.}}

\maketitle

\begin{abstract}
In this paper, both dynamic mean-variance portfolio selection problems and dynamic variance hedging problems are discussed under non-Markovian framework. Explicit closed-loop equilibrium strategies of these problems are respectively obtained via a unified approach for the first time. Several new interesting facts arise in mean-variance problems with constant risk aversion. For example, it is shown that equilibrium strategies are still allowed to rely on initial wealth as long as risk-free return rate is random. In addition, the closed-loop equilibrium strategy and open-loop equilibrium strategy have the same connection with initial wealth in non-Markovian setting, and they happen to equal to each other if only risk-free return rate is deterministic.
%
%
 
%
%
%

%
\end{abstract}

\ms


\ms



\section{Introduction}

The well-known Markowitz's mean-variance portfolio selection has laid down the foundation of modern
investment portfolio theory. If one considers this financial problem in multi-period setting, one will encounter the so-called time-inconsistency issue. That is to say, the
strategy at this moment may not keep optimality at next moment, which implies people have to change the strategy in a silly way all the time to keep optimality.

Besides optimality, the time consistency of policies is also a basic requirement
for rational decision making in many situations (\cite{Strotz-1956}). Recently, the notion of time consistent equilibrium investment strategies of dynamic mean-variance portfolio selection problems have attracted much attention in the literature. Basically speaking, there are two types of equilibrium strategies along this line: open-loop equilibrium (investment) strategies and closed-loop equilibrium (investment) strategies.
We first give a revisit to the investigations on the former notion. It was introduced and investigated carefully in \cite{HJZ 2012}, \cite{HJZ 2015} when partial involved coefficients are random. Later, \cite{Wei-Wang-2017} extended it into the general asset-liability management problem with full random coefficients. More related topics can also be found in \cite{Alia-Chighoub-Sohail-2016}, \cite{Wang-Wu-2016}, and the references therein. As to the closed-loop equilibrium strategies, they were firstly developed in \cite{BC 2010} under the Markovian framework. To avoid the wealth-independent property of equilibrium strategies in \cite{BC 2010}, the authors in \cite{BMZ 2012} introduced one kind of state dependent risk aversion, and derived the equilibrium strategies via the extended HJB equations idea developed in \cite{Bjork-Murgoci 2014}. We refer to \cite{C-2013}, \cite{Wei-Yong-Yu-2017}, \cite{Huang-Li-Wang-2017}, \cite{Yong-Tams-2017}, \cite{Zeng-Li-2011}, etc., for more related topics.

Mean-variance portfolio selection problems also have inherent connections with variance hedging problems. In financial market, people (or hedgers) always face the risk of non-tradable or a contingent claim at some future time. To hedge this uncertainty, they can do continuous trading in non-risky asset (e.g. bond) and risky asset (e.g. stock). In the literature, minimum variance criterion is widely employed to measure the quality of the hedging, see e.g. \cite{Duffie-Richardson-1990}, \cite{Schweizer-1994}. In the dynamic setting, similar as mean-variance problems, the optimal hedging strategies become time inconsistent as well. In \cite{Basak-Chabakauri-2012}, the author discussed the time consistent equilibrium hedging strategy in Markovian setting.

As far as we know, there is no literature on discussing closed-loop equilibrium strategies of dynamic mean-variance optimization problems in non-Markovian setting. As to the analogue study of dynamic variance hedging problems, it is also open to our best. Motivated by these facts, in this paper we make an attempt to fill these gaps. Compared with the references aforementioned, the approach developed here is advantageous in two aspects. In the first place, it is quite challenging to extend the HJB equation idea in \cite{Bjork-Murgoci 2014}, \cite{BMZ 2012} or the multi-player differential games method in \cite{Yong-Tams-2017} into the non-Markovian setting, not to mention the complicated convergence arguments in these papers. Inspired by \cite{Huang-Li-Wang-2017}, \cite{Wang-AMO-2018}, here we employ another approach and obtain the desired Riccati system without the delicate convergence arguments. In the second place, the introduced method enables us to give a unified treatment on mean-variance portfolio selection problems and variance hedging problem with random coefficients. In addition, the obtained results, which appear for the first time, cover the analogue version in \cite{BC 2010}, \cite{Huang-Li-Wang-2017}, \cite{Bjork-Murgoci 2014}, \cite{BMZ 2012}.

There are also several interesting facts which are revealed here for the first time.

(1) It is well-known the equilibrium strategies of mean-variance problems with deterministic coefficients keep unchanged, no matter what the initial wealth is. However, in our non-Markovian scenario, it is a feedback form of equilibrium wealth process, and actually relies on initial wealth.

(2) In the existing papers with constant risk aversion (e.g., \cite{BC 2010}, \cite{BMZ 2012}, \cite{Huang-Li-Wang-2017}, \cite{Zeng-Li-2011}), the closed-loop equilibrium strategies of mean-variance problem do not depend on initial wealth, which inspires the introducing of state dependent risk aversion (\cite{BMZ 2012}). However, thanks to the randomness of coefficients, we obtain the state dependent equilibrium strategies with merely constant risk aversion.

(3) For the mean-variance problems, we found that both closed-loop equilibrium strategies and open-loop equilibrium strategies have the same way in depending on initial wealth, even when all the coefficients are random. If the risk-free return rate is deterministic, these two equilibrium strategies equal to each other if they exist.

(4) In non-Markovian setting, both the equilibrium investment strategy of mean-variance problem and the equilibrium hedging strategy of variance hedging problem have the same relations with initial wealth.

The rest of this paper is organized as follows. In Section 2, some notations, spaces are introduced and the two financial problems are formulated in detail. Section 3 includes three parts. The first part is aim to study an auxiliary time inconsistent optimal control problem, and provide proper sufficient conditions for closed-loop equilibrium operators that are defined in Section 2. The second part and third part are devoted to treating the mean-variance problem and variance hedging problem respectively. In Section 4, some concluding remarks are present.

\section{Preliminary notations and model formulation}

Through this paper, let $(\Omega, \cF,
P, \{\cF_t\}_{t\geq0})$ be filtered complete probability space, $\{W(t),
t\geq 0\}$ be $\dbF:=\{\mathcal{F}_t\}_{t\geq 0}$-adapted one-dimensional Brownian
motion.

For $ n, p\in\dbN,$ $0\le s<t\le T$, we define
$$\ba{ll}
\ns\ds L^2_{\cF_t}(\Omega;\dbR^n):=\Big\{X:\Omega\to \dbR^n\bigm| X\hb{
is $\cF_t$ measurable,}\ \dbE|X|^2<\infty\Big\},\\
\ns\ds L^2_{\dbF}\big(\O;C([s,t];\dbR^n)\big):=\Big\{X:[s,t]\times\O\to
\dbR^n\bigm|X(\cd)\hb{ is measurable, $\dbF$-adapted, }\\
\ns\ds\qq \qq\qq\qq\qq \hb{has continuous path,}\ \dbE\(\sup_{r\in[s,t]}|X(r)|^2\)<\infty\Big\},\\
\ns\ds L^\infty_{\dbF}(\Omega;C([s,t];\dbR^n)):=\Big\{X:[s,t]\times\Omega\to \dbR^n\bigm|X(\cd)\hb{
is  measurable and}\\
\ns\ds\qq \qq\qq\qq\qq \hb{$\dbF$-adapted, has continuous paths}, \ ~\esssup_{\omega\in\Omega}\sup_{r\in[s,t]}|X(r)| <\infty\Big\},\\
\ns\ds L^p_{\dbF}(\Omega;L^2(s,t;\dbR^n)):=\Big\{X:[s,t]\times\O\to
\dbR^n\bigm|X(\cd)\hb{ is measurable, $\dbF$-adapted,}\\
\ns\ds\qq \qq\qq\qq\qq \ \dbE\(\int_{s}^t|X(r)|^2dr\)^{\frac p 2}<\infty\Big\},\\
\ns\ds
 L^p_{\dbF}(\Omega;L^1(s,t;\dbR^n)):=\Big\{X:[s,t]\times\O\to
\dbR^n\bigm|X(\cd)\hb{ is measurable, $\dbF$-adapted,}\\
\ns\ds\qq \qq\qq\qq\qq \ \dbE\(\int_{s}^t|X(r)|dr\)^{  p }<\infty\Big\}.
\ea$$
In particular, $L^2_{\dbF}(s,t;\dbR^n)=L^2_{\dbF}(\Omega;L^2(s,t;\dbR^n))$.

We consider a financial market where two assets are traded continuously on $[0,T]$. Suppose the price of bond evolves as
\begin{eqnarray*}
\left\{\begin{array}{rl}
dS_0(s) & \!\!\!= r(s)S_0(s)ds,\;\; s \in [0, T],\\
 S_0(0) & \!\!\!= s_0>0,
\end{array}\right.
\end{eqnarray*}
and the risky asset is described by
\begin{eqnarray*}
\left\{\begin{array}{rl} dS(s) & \!\!\!=
S(s)\Big\{b(s)ds+ \si (s)dW(s)\Big\},\;\;
s \in [0, T],\\
 S(0) & \!\!\!= s_1>0.\;\;
\end{array}\right.
\end{eqnarray*}
Here $r >0$ is the risk-free return rate, $b $ is the expected return rate of risky asset, $\si $ is the corresponding volatility rate.

\ms

(H0) Suppose $r,$ $b $, $\sigma$ are bounded and $\dbF$-adapted processes, and there exists constant $\d>0$ such that $|\sigma|^2 \geq \delta.$
%

\ms

Given initial capital $x>0$, $\beta:=b-r$, $\theta :=\beta \si^{-1}$, for $s\in[0,T]$, the investor's wealth $X(s)$ satisfies
\begin{equation}\label{wealth-equation}
\!\!\!\!\!\!\left\{\begin{array}{rl}
\!\!\!dX(s)  & \!\!\!=\big[r(s) X(s)+\beta(s)  u(s) \big] ds+u(s)\si(s)  dW(s),  \\
 \!\!\!X(0) & \!\!\!=x,
\end{array}\right.
\end{equation}
where $u $ is the capital invested in the risky asset.

At time $t$, the objective of a mean-variance portfolio selection problem is to choose an investment strategy to minimize
\begin{eqnarray}\label{cost-functional}
J(u(\cd);t,X(t))=\hb{Var}_{t}\big[X(T)\big]-\gamma\dbE_{t}\big[X(T)\big],
\end{eqnarray}
where $\dbE_{t}[\cdot]:=\dbE[\cdot|\mathcal{F}_t],$
$\gamma$ is constant risk aversion.

Inspired by the existing papers (\cite{Huang-Li-Wang-2017}, \cite{Wang-AMO-2018}, \cite{Wei-Wang-2017}), we introduce the following definition in non-Markovian setting.

To this end, given $t\in[0,T]$, $\e>0$, proper $(\Th^*,\f^*)$ and bounded $v\in L^2_{\cF_t}(\Omega;\dbR)$, let
\bel{Definitions-of-equilibrium-perturbation}
\ba{ll}
\ns\ds u^{v,\e}:=\Th^* X^{v,\e}+\f^*+ vI_{[t,t+\e]}, \ \ u^{*}:=\Th^* X^{*}+\f^*,
\ea
\ee
where for $s\in[0,T]$, $(X^{v,\e}(s),X^*(s))$ are described as
\bel{Equilibrium-state-closed-loop}\left\{\ba{ll}
\ns\ds dX^*(s)=\Big[\big[r(s)+\b(s)\Th^*(s)\big]X^*(s)+\b(s)\f^*(s)\Big]ds\\
\ns\ds\qq\qq\q +\si(s)(\Th^*(s)X^*(s)+\f^*(s))
dW(s),\\
\ns\ds dX^{v,\e}(s)=\Big[\big[r(s)+\b(s)\Th^*(s)\big]X^{v,\e}(s)+\b(s)\f^*(s)+\b(s) vI_{[t,t+\e]}(s)\Big]ds\\
\ns\ds\qq\qq\q +\si(s)\big[ \Th^*(s)X^{v,\e}(s)+\f^*(s) +vI_{[t,t+\e]}(s)\big]
dW(s),\\
\ns\ds X^*(0)=x, \ \ X^{v,\e}(0)=x.
\ea\right.
\ee

\bde\label{Definition-closed-strategy}
A pair of processes $(\Th^*(\cd),\f^*(\cd))\in
L^p_{\dbF}(\Omega;L^2(0,T;\dbR))\times L^2_{\dbF}(0,T;\dbR)$ is called a closed-loop equilibrium operator if for any $x\in\dbR$, $u^*, u^{v,\e}\in L^2_{\dbF}(0,T;\dbR)$, and
\bel{optimal-open}\lim_{\overline{\e\to0}}
{J(u^{v,\e}(\cd);t,X^*(t))-J\big(u^*(\cd)\big|_{[t,T]};t,X^*(t)\big)
\over\e}\ge0.
\ee
Above $u^*$, $X^*$ is called closed-loop equilibrium investment strategy, closed-loop equilibrium wealth process, respectively.
\ede

\ms

Next we discuss the formulation of variance hedging problems. An agent (or hedger) is committed to hold a non-tradable asset with payoff $\xi$ at time $T$.
Here the
asset may be interpreted as a derivative security or a real option, and $\dbE|\xi|^k<\infty$, $k>2$.
The corresponding risk can be hedged by continuous trading in bond and risky asset aforementioned. The hedging policy $\pi$ is the dollar amount invested in the stock. Hence the hedger's tradable wealth $X$ is described by (\ref{wealth-equation}) with $u$ replaced by $\pi$.
The aim is to find a proper hedging policy to minimize
\bel{}\ba{ll}
\ns\ds \h J(\pi(\cd);t,X(t)):=\hb{Var}_{t}\big[\xi-X(T)\big].
\ea\ee

Let $\l(\cd):=\dbE\big[\xi\big|\cF_\cd\big]$, $Y(\cd):=X(\cd)-\l(\cd)$, we have
\bel{Variance-hedging-new-state}\left\{\ba{ll}
\ns\ds  dY(s)=\big[r(s)Y(s)+r(s)\l(s)+\b(s)\pi(s)\big]ds\\
\ns\ds\qq\qq +\big[\pi(s)\si(s)-\zeta(s)\big]dW(s),\ \ s\in[0,T],\\
\ns\ds d\l(s)=\zeta(s)dW(s),\ \ s\in[0,T],\\
\ns\ds Y(0)=x-\dbE\xi,\ \ \l(T)=\xi.
\ea\right.\ee
Notice that both the diffusion term and drift term of $Y(\cd)$ include the nonhomogeneous terms.
The aim of the hedger is to minimize $\h J(\pi(\cd);t,X(t))=\dbE_t \big|Y(T)\big|^2-\big|\dbE_t Y(T)\big|^2$.

Similar as Definition \ref{Definition-closed-strategy},
 we can define the closed-loop equilibrium operator $(\Th^*,\f^*)$, closed-loop equilibrium hedging policy $\pi^*$ as well. We omit it for simplicity.

In the sequel, $K$ is a generic constant which varies in different context.
\ms

%

\section{Equilibrium strategies in mean-variance problems and variance hedging problems}

\subsection{An auxiliary time inconsistent optimal control problem}

To give a unified treatment of above financial problems, we study a slightly general optimal control problem associated with (\ref{cost-functional}) and
\bel{General-state-equation}\left\{\ba{ll}
\ns\ds
dX(s)=\big[r(s)  X(s) +\b (s)u(s)+l(s) \big]ds   +\big[ \si (s) u(s) +h(s)\big]dW(s), \ \ s\in[0,T],\\
\ns\ds X(0)=x.
\ea\right.
\ee
Similar as Definition \ref{Definition-closed-strategy}, we define the corresponding closed-loop equilibrium operator as well.

In the following, we assume that

\ms

(H1) For any $p>2$, $h\in L^p_{\dbF}(\Omega;L^2(0,T;\dbR))$, $l\in L^p_{\dbF}(\Omega;L^1(0,T;\dbR))$.

\ms

For proper $(\Th^*,\f^*)$, constant $\e>0$, $\cF_t$-measurable bounded random variable $v$, we define
\bel{Definition-of-u-0-v-episilon}
\ba{ll}
\ns\ds u^*:= \Th^* X^*+\f^*,\ \ u^{v,\e}_0:=\Th^* X^{v,\e}_0+\f^*+vI_{[t,t+\e]},
\ea
\ee
where for $s\in[0,T],$ $X^{*}(s)$, $X^{v,\e}_0(s)$ are described as,
\bel{State-equation-Th-1-Th-2}\left\{\2n\ba{ll}
\ns\ds
dX^*(s)=\big[r(s)  X^*(s) +\b (s)u^*(s)+l(s) \big]ds   +\big[ \si (s) u^*(s) +h(s)\big]dW(s), \ \ \\
\ns\ds
dX^{v,\e}_0(s)=\big[ r(s) X^{v,\e}_0(s)+\b(s)u^{v,\e}_0(s)+l(s)\big]ds +  \big[\si(s) u^{v,\e}_0 (s)+h(s)\big]dW(s),  \\
\ns\ds  X^*(0)=x,\ \ X^{v,\e}_0(0)=x. \ \
\ea\right.\ee

\ms

(H2) Suppose there exists
$$(\Th^*,\f^*)\in L^p_{\dbF}(\Omega;L^2(0,T;\dbR))\times L^2_{\dbF}(0,T;\dbR)$$
such that  for any $x\in\dbR$, $u^*, u^{v,\e}_0\in L^2_{\dbF} (0,T;\dbR)$.

\ms

Under assumption (H2), $X^* $, $X^{v,\e}_0\in L^2_{\dbF}(\Omega;C([0,T];\dbR))$, and the following are well-defined, %
$$J(u^*;t,X^*(t)),\ \ J(u^{v,\e}_0;t,X^*(t)),\ \ t\in[0,T).$$

Given $(\Th^*,\f^*)$ in (H2), we introduce four BSDEs on $[0,T]$ as follows:
\bel{Equations-for-P-i}\left\{\ba{ll}
\ns\ds dP_1^* =-\Big\{2r P_1^* +2(P_1^* \b +\L_1^* \si )\Th^*+|\Th^*|^2\si^2P_1^* \Big\}ds+\L_1^* dW(s),\\
%
\ns\ds dP_2^* =- \big[ r  P_2^* +\Th^*( \b P_2^* +\si\L_2^*) \big]ds+
\L_2^* dW(s),\\
\ns\ds dP_3^* =- \big[ (P_2^* \b +\L_2^*\si)\f^* +P_2 ^*l+\L_2^* h\big]ds+\L_3^*dW(s),\\
\ns\ds dP_4^* =-\big[(r+\Th^*\b)P_4^*+\Th^*\si\L_4^*+( P_1^*\b+\L_1^*\si+
\Th^*\si^2P_1^*)\f^* \\
\ns\ds\qq\q +P_1^* l+\L_1^* h+\Th^* \si P_1^* h\big]ds+
\L_4^* dW(s),\\
\ns\ds P_1^*(T)=2,\ \ P_2^*(T)=1, \ \ P_3^*(T)=0,\ \ P_4^*(T)=-\g.
\ea\right.\ee

We will discuss the solvability of (\ref{Equations-for-P-i}) later. Before that, let us provide a sufficient condition for closed-loop equilibrium operator of optimal control problem associated with (\ref{General-state-equation}), (\ref{cost-functional}).

\ms

\bt\label{Main-theorem-1}
Suppose (H0) holds, and there exist four pairs of processes $(P_i^*,\L_i^*)$ satisfying system (\ref{Equations-for-P-i}) where $P_1^*>\d>0,$ $\d$ is a constant, and
\bel{Definitions-Th-varphi-1}\left\{\ba{ll}
\ns\ds \Th^*:=-\frac{\b(P_1^*-2|P_2^*|^2)+\si(\L_1^*-2\L_2^*P_2^*) }{\si^2 P_1^*},\\
\ns\ds \f^*:=-\frac{\b(P_4^*-2P_2^*P_3^*)+\si(P_1^*h+\L_4^*-2\L_2^*P_3^*)}{\si^2P_1^*}.
\ea\right.\ee
Moreover, (H2) holds, and for any $p>2$,
\bel{Regularity-of-P-i}\left\{\ba{ll}
\ns\ds  (P_i^*,\L_i^*)\in L^{\infty}_{\dbF}(\Omega;C([0,T];\dbR))\times L^{p}_{\dbF}(\Omega;L^2(0,T;\dbR)),\ \ i:=1,2, \\
\ns\ds  (P_j^*,\L_j^*)\in L^{p}_{\dbF}(\Omega;C([0,T];\dbR))\times L^{p}_{\dbF}(\Omega;L^2(0,T;\dbR)),\ \ j:=3,4,\\
\ns\ds\qq \qq \qq \sup_{t\in[\t,T]}\dbE_\t|\L_2^*(t)|^2<\infty,\ \ a.s.\ \  \t\in[0,T].
\ea\right.\ee
Then $(\Th^*,\f^*)$ is a closed-loop equilibrium operator.
\et
\br
As to $\L_2(\cd)$, above pointwise integrability in (\ref{Regularity-of-P-i}) will play an important role next, even though it's stronger than the conventional square integrability. We will verify it later.
\er

\it Proof. \rm \it Step 1. \rm For $t\in[0,T)$, $(\Th^*,\f^*)$ in (\ref{Definitions-Th-varphi-1}), $X^*$ in (\ref{State-equation-Th-1-Th-2}), we introduce
\bel{MF-BSDEs-*-1}\left\{\ba{ll}
\ns\ds dY^*(s,t)=  -\big[(r(s)+\b(s)\Th^*(s))Y ^*(s,t)+\Th^*(s) \si(s) Z^*(s,t)\big]ds\\
\ns\ds\qq\qq\q +Z^*(s,t)dW(s),\ \ s\in[t,T],\\
\ns\ds Y^*(T,t)=2X^*(T)- 2\dbE_tX^*(T)-\g.
\ea\right.\ee
In this step, for $t\in[0,T]$, we prove that the following pair of processes satisfy (\ref{MF-BSDEs-*-1}),
\bel{Representation-for-Y-*}\left\{\ba{ll}
\ns\ds  Y'(\cd,t):=P_1^*X^*- 2P_2^*\dbE_t\big[P_2^*X^*+P_3^*\big]+ P_{4}^*, \  \ \\
\ns\ds Z'(\cd,t):=-2\L_2^*\dbE_t[P_2^*X^*+P_3^*] +(P_1^*\si\Th^* +\L_1^*) X^*+ P_1^*\si\f^*+P_1^* h+\L_4^*.
\ea\right.
\ee

In fact, by It\^{o}'s formula,
$$\left\{\ba{ll}
\ns\ds
d P_1^*X^* =\big[P_1^* (r+\b\Th^*)+\Pi_1+\L_1^* \si\Th^*
\big]X^* ds + (P_1^*[ \si\Th^* X^*+ \si\f^*+h]+\L_1^*X^*) dW(s)\\
\ns\ds\qq\qq\qq
+ \big[(P_1^*\b+\L_1^*\si )\f+P_1^*l+\L_1^*h\big] ds,\\
\ns\ds d\Big[-2P_2^*\dbE_t[ P_2^*X^* ]\Big] =\Big\{-2\Pi_2\dbE_t[
P_2^*X^* ]-2P_2^*\Big\{\dbE_t\big[ (P_2^*r+P_2^*\b\Th^* + \Pi _2+\1n \L_2^* \si\Th^*)X^* \big]
\\
\ns\ds\qq\qq\qq\qq\qq +\dbE_t\big[(P_2^*\b+\L_2^*\si)\f^* +P_2^*l +\L_2^*h \big] \Big\}\Big\}ds-2\L_2^*
\dbE_t[P_2^*X^*]dW(s),\\
\ns\ds d\Big[-2 P_2^*\dbE_t P_3^*\Big] =-2\big[\Pi_2\dbE_t P_3^*+P_2^* \dbE_t \Pi_3 \big]ds-2\L_2^*
\dbE_t P_3^* dW(s),
\ea\right.$$
where $\Pi_i$ $(1\leq i\leq 4)$ denotes the generator of BSDEs in (\ref{Equations-for-P-i}) respectively.
Consequently,
$$\ba{ll}
\ns\ds dY'=d\big[P_1^*X^*-2P_2^*\dbE_t(P_2^*X^*+P_3^*)+P_4^*\big]\\
\ns\ds \qq =\Big\{\big[P_1^* (r+\b\Th^*)+\Pi_1+\L_1^* \si\Th^*
\big]X^*+ (P_1^*\b+\L_1^*\si )\f^* +P_1^* l+\L_1^* h\\
\ns\ds\q\qq -2 \Pi_2\dbE_t[
P_2^*X^*+P_3^*]-2P_2^*\Big\{\dbE_t\big[ (rP_2^*+\b P_2^*\Th^* + \Pi _2+\1n \L_2^* \si \Th^*)X^* \big]
\\
\ns\ds\q\qq+\dbE_t\big[(\b P_2^*+\L_2^*\si)\f^*+P_2^*l+\L_2^*h+\Pi_3\big] \Big\} +\Pi_4\Big\}ds\\
\ns\ds \q\qq + \Big\{-2\L_2^*\dbE_t[P_2^*X^*] + P_1^*[ \si\Th^* X^*+ \si\f^*+h] -2\L_2^*
\dbE_t P_3^*+\L_4^*+\L_1^* X^*\Big\} dW(s).
\ea$$
Putting the definitions of $\Pi_i$ into above equality, we immediately obtain the conclusion.

\ms

\it Step 2. \rm
For $(Y',Z')$ in (\ref{Representation-for-Y-*}), we see that
$$(\dbY'(s),\dbZ'(s)):=(Y'(s,s),Z'(s,s)),\ \ s\in[0,T],$$
are well-defined. In this step, we prove that
\bel{Equality-0}\ba{ll}
\ns\ds \lim_{\e\rightarrow0}\Big[\frac 1 \e \dbE_t\int_t^{t+\e}\big[\b(s) Y'(s,t)+\si(s) Z'(s,t)\big]ds\Big]\\
\ns\ds = \lim_{\e\rightarrow0}\Big[\frac 1 \e \dbE_t\int_t^{t+\e}\big[\b(s) \dbY'(s)+\si(s) \dbZ'(s)\big]ds\Big]=0.
\ea
\ee

Notice that the second equality follows from (\ref{Definitions-Th-varphi-1}). We focus on the first one.

By the definitions of $(Y',Z')$ in (\ref{Representation-for-Y-*}),
$$\ba{ll}
\ns\ds \b  Y'(\cd,t)+\si Z'(\cd,t)\\
\ns\ds =\b\big[P_1^* X -2 P_2^* \dbE_t\big[P_2 ^*X^* +P_3 ^*\big]+ P_{4}^* \big]+
\si\big[-2\L_2^*\dbE_t[P_2^*X^*+P_3^*] \\
\ns\ds\q  +(P_1^*\si\Th^* +\L_1^*) X^*+ P_1^*\si\f^*+P_1^* h+\L_4^*\big]\\
\ns\ds=\big[\b P_1^*+\si(P_1^*\si\Th^*+\L_1^*)\big]X^*+\si^2 P_1^*\f^*+\si P_1^* h+\b P_4^*+\si\L_4^*\\
\ns\ds\q -2\big[\b P_2^*+\si\L_2^*\big] \dbE_t(P_2^* X^*+P_3^*).
\ea$$
In particular, one has
\bel{Diagonal-value-summation}\ba{ll}
\ns\ds  \b \dbY' +\si \dbZ'   =\Big[\b(P_1^*-2|P_2^*|^2)+\si(-2\L_2^*P_2^*+P_1^*\si\Th^*+\L_1^*)\Big]X^*\\
\ns\ds \qq \qq\qq +\b(-2P_2^*P_3^*+P_4^*)+\si(-2\L_2^*P_3^*+P_1^*\si\f^*+P_1^*h+\L_4^*).
\ea\ee
Therefore,
\bel{Relation-dbY-Y'-0}\ba{ll}
\ns\ds \dbE_t\big[\b Y'(\cd,t)+\si Z'(\cd,t)\big]\\
\ns\ds =
 \dbE_t\big[\b \dbY' +\si \dbZ' \big]
-2\dbE_t\Big\{\big[\dbE_t( P_2^*X^*+P_3^*)-P_2^*X^*-P_3^*\big] (\b P_2^*+\si\L_2^*)\Big\}.
\ea\ee
In terms of (\ref{Regularity-of-P-i}) and conditional dominated convergence theorem, we conclude that
$$\ba{ll}
\ns\ds   \lim_{\e\rightarrow0}\dbE_t\sup_{s\in[t,t+\e]}\big|\dbE_t(P_2^*(s)X^*(s)+P_3^*(s))-(P_2^*(s)X^*(s)+P_3^*(s))\big|
=0. \ \ a.s.
\ea$$
Then as $\e\rightarrow0$, one has
\bel{Result-convergence-L-2}\ba{ll}
\ns\ds \frac{1}{\e}
\dbE_t\int_t^{t+\e}\big|\L_2^*(s)\big[\dbE_t(P_2^*(s)X^*(s)+P_3^*(s))-(P_2^*(s)X^*(s)+P_3^*(s))\big]\big|ds\\
\ns\ds\leq \sup_{s\in [t,T]}\big[\dbE_t |\L_2^*(s)|^2 \big]^{\frac 1 2}\Big[\dbE_t\sup_{s\in[t,t+\e]}\big|\dbE_t(P_2^*(s)X^*(s)+P_3^*(s))-(P_2^*(s)X^*(s)+P_3^*(s))\big|^2\Big]
^{\frac 1 2}\rightarrow0.
\ea\ee
Similarly, by the imposed regularity of $P_2^*$,
\bel{Result-convergence-P-2}\ba{ll}
\ns\ds \frac{1}{\e}
\dbE_t\int_t^{t+\e}\big|P_2^*(s)\big[\dbE_t(P_2^*(s)X^*(s)+P_3^*(s))-(P_2^*(s)X^*(s)+P_3^*(s))\big]\big|ds\rightarrow0,\ \ \e\rightarrow0.
\ea\ee
Putting (\ref{Result-convergence-L-2}) and (\ref{Result-convergence-P-2}) back into (\ref{Relation-dbY-Y'-0}), we get the desired conclusion.

 \ms

 \it Step 3. \rm For $s\in[0,T]$, we define $X_1^{v,\e}(s):=X^{v,\e}_0(s)-X^*(s)$ that satisfies
\bel{Equation-of-s-X-1-epislon}\left\{\2n\ba{ll}
\ns\ds
dX_1^{v,\e}(s) =\big[ (r (s)+\b(s)\Th^*(s))  X_1^{v,\e}(s)   +\b (s) vI_{[t,t+\e]}(s)  \big]ds\\
 \ns\ds\qq \qq \q+ \big[\si(s)\Th^*(s) X_1^{v,\e}(s)+ \si (s) vI_{[t,t+\e]} (s) \big]dW(s), \\
\ns\ds X^{v,\e}_1(0)=0.\ea\right.\ee
From the definition of $u^{v,\e}_0$ in (\ref{Definition-of-u-0-v-episilon}), it is a direct calculation that
\bel{Inequality-for-Y-Z}\ba{ll}
\ns\ds J(u^{v,\e}_0(\cd);t,X^*(t))-J(u(\cd);t,X^*(t)) \\
\ns\ds = \dbE_t\Big\{\big[2X^*(T)- 2\dbE_tX^*(T)- \g  \big] X_1^{v,\e}(T)\Big\}\\
\ns\ds\qq +\dbE_t \Big\{( X_1^{v,\e}(T)- \dbE_tX_1^{v,\e}(T))X_1^{v,\e}(T)\Big\}\\
\ns\ds\geq  \dbE_t\Big\{\big[2X^*(T)- 2\dbE_tX^*(T)- \g \big] X_1^{v,\e}(T)\Big\}.
\ea\ee
The following result is implied by It\^{o}'s formula to $Y'(\cd,t)X_1^{v,\e}(\cd)$,
\bel{Ito-result-1} \ba{ll}
\ns\ds \dbE_t\Big\{\big[2X^*(T)- 2\dbE_tX^*(T)-\g \big] X_1^{v,\e}(T)\Big\}=\dbE_t \int_t^{t+\e}(\b(s)Y'(s,t)+\si(s)Z'(s,t))ds \cd  v.
\ea
\ee
To sum up, from (\ref{Equality-0}), (\ref{Inequality-for-Y-Z}), (\ref{Ito-result-1}), we conclude that
$$\ba{ll}
\ns\ds \lim_{\e\rightarrow0}\frac{J(u^\e(\cd);t,X^*(t))-J(u^*(\cd);t,X^*(t))}{\e} \\
\ns\ds\geq \lim_{\e\rightarrow0}\Big[\frac 1 \e \dbE_t\int_t^{t+\e}\big[ \b(s)Y'(s,t)+\si(s) Z'(s,t) \big] v ds \\
\ns\ds =\lim_{\e\rightarrow0}\Big[\frac 1 \e \dbE_t\int_t^{t+\e}\big[\b(s)  \dbY'(s) +\si(s) \dbZ'(s) \big]ds\Big] v=0.
\ea$$
\endpf

\ms

The following result ensures the requirements imposed in Theorem \ref{Main-theorem-1}. 

\ms

(H3) For $s\in[0,T]$, $r(s)$, $\frac{\b(s)}{\si(s)}$ are Malliavin differentiable, and there exists constant $K>0$ such that $\Big[\big|D_{\nu}r(s)\big|+\big|D_{\nu}\big[\frac{\b(s)}{\si(s)}\big]\big|\Big]\leq K,$ $\nu, s\in[0,T]$.

\bt\label{Theorem-main-2}
Suppose (H0), (H1), (H3) hold. Then there exists four pairs of $(P_i^*,\L_i^*)$ satisfying system (\ref{Equations-for-P-i}) and condition (\ref{Regularity-of-P-i}) with
\bel{True-Th-varphi}\left\{\ba{ll}
\ns\ds P_1^*=2|P_2^*|^2,\ \ \L_1^*=4P_2^*\L_2^*,\ \ \Th^*:=-\frac{\L_2^*}{\si P_2^*},\ \ \\
\ns\ds \f^*:=-\frac{2P_2^*\L_3^*}{\si P_1^*}-\frac{(\b\Phi+\si\Psi+\si P_1^*h)}{\si^2 P_1^*},\\
\ns\ds \f^*\in L^p_{\dbF}(\Omega;L^2(0,T;\dbR)),\ \ p>2.
\ea\right.\ee
Here
\bel{Equation-of-Phi-Psi}\left\{\ba{ll}
\ns\ds d\Phi=-\Big[(r+\b\Th^*)\Phi+\Th^*\si\Psi \Big]ds+\Psi dW(s),\\
\ns\ds \Phi(T)=-\g.
\ea\right.
\ee
Moreover, $(\Th^*,\f^*)$ can be rewritten as the forms in (\ref{Definitions-Th-varphi-1}), and assumption (H2) is fulfilled.
\et

\it Proof. \rm
\it Step 1. \rm We prove that system (\ref{Equations-for-P-i}) is solvable associated with (\ref{Definitions-Th-varphi-1}) and (\ref{Regularity-of-P-i}).

At first, we consider
\bel{One-construction-M-N-1}\left\{\ba{ll}
\ns\ds dM=\big[r M+\frac{\b}{\si}N\big]ds+ N dW(s),\ \ s\in[0,T],\\
\ns\ds M(T)=-\sqrt{2}.
\ea\right.\ee
It is easy to see $M$ is bounded and $\int_0^{\cd}N(s)dW(s)$ is a BMO-martingale.

By defining $P_1^*:=\frac{4}{M^2},$ $\L_1^*:=-\frac{8N}{M^3}$, we have
\bel{Equation-for-P-L-1}\left\{\ba{ll}
\ns\ds  dP_1^*=-\Big\{2r P_1^* -
\frac{\b\L_1^*}{\si} -\frac{3|\L_1^*|^2}{4P_1^*} \Big\}ds+\L_1 ^*dW(s),\ \ s\in[0,T],\\
\ns\ds P_1^*(T)=2.
\ea\right.\ee
Since $M$ is bounded, $\int_0^{\cd}\L_1(s)dW(s)$ is BMO-martingale.

 We define
\bel{several-definitions}
\ba{ll}
\ns\ds P_2^*:=\sqrt{\frac{P_1^*}{2} },\ \ \L_2^*:=\frac{\L_1^*}{4P_2^*},\ \ \Th^*:=-\frac{2\si\L_2^*P_2^*}{\si^2 P_1^*}.
\ea
\ee
Therefore, we obtain the first expression in (\ref{Definitions-Th-varphi-1}), and the first two results in (\ref{True-Th-varphi}).
Moreover, $\int_0^{\cd}\L_2^*(s)dW(s)$, $\int_0^{\cd}\Th^*(s)dW(s)$ are BMO martingales and
$$\ba{ll}
-  \frac{  \b\L_1^*}{\si} -\frac{3|\L_1^*|^2}{4P_1^*} =-\frac{4(P_1^* \b +\L_1^* \si ) \si\L_2^*P_2^*}{\si^2P_1^*} +\frac{\L_2^*P_2^*\L_1^*}{P_1^*}=2(P_1^* \b +\L_1^* \si )\Th^*+|\Th^*|^2\si^2P_1^*.
\ea
$$
As a result, we can rewrite (\ref{Equation-for-P-L-1}) as the first equation in (\ref{Equations-for-P-i}).

To obtain the case of $(P_2^*,\L_2^*)$, we first use It\^{o}'s formula as follows,
\bel{Equation-for-P-3-0}\ba{ll}
\ns\ds  d\big[\frac{P_1^*}{2}\big]^{\frac 1 2} =\Big[ -r(\frac {P_1^*}{2})^{\frac 1 2}+2^{-\frac 3 2}|P_1^*|^{-\frac 1 2}\frac{\b}{\si}\L_1^*+2^{-\frac 5 2}|P_1^*|^{-\frac 3 2}|\L_1^*|^2\Big]ds+2^{-\frac 3 2}|P_1^*|^{-\frac 1 2}\L_1^*dW(s).
\ea\ee
We observe that
$$
(\b P_2^*+\si\L_2^*)\Th^*=-\big[\b (\frac{P_1^*}{2})^{\frac 1 2}+\frac{\si\L_1^*}{4P_2^*}\big]\frac{\L_1^*}{2\si P_1^*}=-\Big[2^{-\frac 3 2}|P_1^*|^{-\frac 1 2}\frac{\b}{\si}\L_1^*+2^{-\frac 5 2}|P_1^*|^{-\frac 3 2}|\L_1^*|^2\Big].
$$
Recalling $P_2^*(T)=1$, we can rewrite (\ref{Equation-for-P-3-0}) into the the second equation in (\ref{Equations-for-P-i}).

The regularity of $(P_i^*,\L_i^*)$, $i=1,2$ in (\ref{Regularity-of-P-i}) is obvious.

We continue to investigate $(P_3^*,\L_3^*)$.

To begin with, we look at BSDE (\ref{Equation-of-Phi-Psi}).
Since $\int_0^{\cd}\Th^*(s)dW(s)$ is BMO-martingale, for any $p>2$, by Theorem 10 in \cite{Briand-Confortola-2008}, one has
$$
\ba{ll}
\ns\ds (\Phi,\Psi)\in L^p_{\dbF}(\Omega;C([0,T];\dbR))\times  L^p_{\dbF}(\Omega;L^2(0,T;\dbR)).
\ea
$$
Given this pair of $(\Phi,\Psi)$, let us consider
\bel{Equation-of-P-4}\left\{\ba{ll}
\ns\ds dP_3^*=\Big[\frac{(P_2^*\b+\L_2^*\si)2P_2^*}{\si P_1^*}\L_3^*+\Pi
\Big]ds+\L_3^*dW(s),\ \ s\in[0,T],\\
\ns\ds P_3^*(T)=0,
\ea\right.
\ee
where
$$\ba{ll}
\ns\ds \Pi:=\frac{(P_2^*\b+\L_2^* \si)(\b\Phi+\si\Psi+\si P_1^*h)}{\si^2 P_1^*}-P_2^*l-\L_2^*h.
\ea
$$
If we define
\bel{Second-definition-f}\ba{ll}
\ns\ds \f^*:=-\frac{2P_2^*\L_3^*}{\si P_1^*}-\frac{(\b\Phi+\si\Psi+\si P_1^*h)}{\si^2 P_1^*},
\ea
\ee
we transform (\ref{Equation-of-P-4}) into the fourth one in (\ref{Equations-for-P-i}).

By the integrability of $(P_2^*,\L_2^*)$, $(\Phi,\Psi)$, $(h,l)$, for any $p>2$, we have $\Pi\in L^p_{\dbF}(\Omega;L^1(0,T;\dbR))$. As a result, by Theorem 10 of \cite{Briand-Confortola-2008}, we obtain the regularity of $(P_3^*,\L^*_3)$ in (\ref{Regularity-of-P-i}).

The integrability of $\f^*$ in (\ref{True-Th-varphi}) is easy to get.

Now we turn to look at $(P_4^*,\L_4^*)$. Recall the second and third equations of (\ref{Equations-for-P-i}), we have the following by It\^{o}'s formula,  
$$\ba{ll}
\ns\ds d(-2P_2^*P_3^*)=-\Big[-2P_2^*(P_2^*\b+\L_2^*\si)\f^*-2(r+\b\Th^*)P_2^*P_3^*-2P_3^*\Th^* \si\L_2^*-2P_2^*(P_2^*l+\L_2^* h)\\
\ns\ds\qq\qq +2\L_2^*\L_3^*\Big]ds-(2P_2^*\L_3^*+2P_3^*\L_2^*)dW(s).
\ea
$$
Let
\bel{Definition-of-P-5-Lambda}\ba{ll}
\ns\ds
P_4^*:=\Phi+2P_2^*P_3^*,\ \L_4^*:=\Psi+2\L_2^*P_3^*+2\L_3^*P_2^*.
\ea\ee
Therefore, it is a direct calculation that
\bel{Equation-of-P-5}\left\{\ba{ll}
\ns\ds dP_4^* =-\big[(r+\Th^*\b)P_4^*+\Th^*\si\L_4^*-2\Th^*\si P_2^*\L_3^*+2P_2^*( P_2^*\b+\L_2^*\si)\f^* \\
\ns\ds\qq\q +2P_2^*(P_2^* l+\L_2^* h)-2\L_2^*\L_3^*\big]ds+\L_4^* dW(s),\\
\ns\ds P_4^*(T)=-\g.
\ea\right.
\ee
We observe that above (\ref{several-definitions}) yields the following results,
$$\Th^*=-\frac{\L_2^*}{\si P_2^*},\ \ \Th^*\si P_2^*\L_3^*+\L_2^*\L_3^*=0,\ \ \L_1^*+\Th^*\si P_1^*=2P_2^*\L_2^*.$$
Consequently,
$$\ba{ll}
\ns\ds -2 \Th^*\si P_2^*\L_3^*-2\L_2^*\L_3^* +2P_2^*( P_2^*\b+\L_2^*\si)\f^* +2P_2^*(P_2^* l+\L_2^* h)\\
\ns\ds =P_1^*\b\f^*+(\L_1^*\si+\Th^*\si^2 P_1^*)\si\f^*+P_1^*l+(\L_1^*+\Th^*\si P_1^*)h.
\ea$$
Plugging it back into (\ref{Equation-of-P-5}), we obtain the fourth equation in (\ref{Equations-for-P-i}).

The regularity of $(P_4^*,\L_4^*)$ in (\ref{Regularity-of-P-i}) is easy to see.

Eventually, by (\ref{Definition-of-P-5-Lambda}), we can rewrite (\ref{Second-definition-f}) as the second form in (\ref{Definitions-Th-varphi-1}).

\ms

\it Step 2. \rm In this step, we prove the integrability of $\L_2^*$.

Given $\t\in[0,T],$ $\nu\in[\t,T]$, by (H3) and Proposition 5.3 in \cite{El-Karoui-Peng-1997}, the Malliavin derivatives $(D_\nu M(s),D_\nu N(s))$ exist, $s\in[0,T]$, and a version is given by
$$\left\{\ba{ll}
\ns\ds (D_\nu M(t),D_\nu N(t))=(0,0),\ \ a.s. \ a.e., \ \ t\in[0,\nu),\\
\ns\ds D_\nu M(t)=-\int_t^T\Big[D_\nu r(s) \cd M(s)
+D_\nu\big[\frac{\b(s)}{\si(s)}\big] \cd N(s)+r(s) D_\nu M(s)\\
\ns\ds\qq\qq\q +\frac{\b(s)}{\si(s)} D_\nu N(s)\Big]ds-\int_t^TD_\nu N(s)dW(s),\ \ t\in[\nu,T].
\ea\right.$$
By classical estimate of BSDEs, for $\nu\in[\t,T]$,
$$\ba{ll}
\ns\ds \dbE_\t\sup_{t\in[\nu,T]}|D_\nu M(t)|^2+\dbE_\t\int_\nu^T |D_\nu N(s)|^2ds\\
\ns\ds\leq K\dbE_\t\Big[\int_\nu^T\big[ |D_\nu r(s)||M(s)|+|D_\nu\big[\frac{\b(s)}{\si(s)}\big]|
 |N(s)|\big]ds\Big]^2.\ \ a.s.
\ea$$
According to (H3), for $\t\in[0,T]$, we arrive at
$$\ba{ll}
\ns\ds \sup_{\nu\in[\t,T]}\dbE_\t|D_\nu M(\nu)|^2\leq K\sup_{\nu\in[\t,T]}\dbE_\t\Big[\int_\nu^T\big[
|M(s)|+  |N(s)|\big]ds\Big]^2<\infty.
\ea$$
Therefore,
$$\sup\limits_{t\in[\t,T]}\dbE_\t|N(t)|^2=\sup\limits_{t\in[\t,T]}\dbE_\t|D_{t}M(t)|^2<\infty. \ \  a.s.$$
Our conclusion is followed by the definition of $\L_1^*$ and notations in (\ref{several-definitions}).

\ms

\it Step 3. \rm We verify the assumptions in (H2).

At first, let us look at the case of $(u^*,X^*)$. For $s\in[0,T]$, recall that
$$\left\{\ba{ll}
\ns\ds
 dX^*(s) =\big[[r+\b\Theta^*]X^* +\b \f^*+l \big]ds  +\left[\si\Theta^* X^*+ \si\f^*+h \right]dW(s),\\
\ns\ds X^*(0)=x.
\ea\right.
$$
The unique solution to the above linear SDE is given by
$$\ba{ll}
\ns\ds
X^*(t)= \Phi(t)x+\Phi(t)\int_{0}^{t}\Phi^{-1}\big[\f^*(\b-\si^2\Th^*)+l-\si\Th^* h\big]ds\\
\ns\ds\qq\qq+\Phi(t)\int_0^t\Phi^{-1}(\si\f^*+h)dW(s),
\ea$$
where $\Phi(\cdot)$ satisfies
\[
\begin{cases}
d\Phi(s)=\Phi( r+\b\Th^*)ds+ \si\Th^*\Phi dW(s),\q s\in[0,T],\\
\Phi(0)=1.
\end{cases}
\]
It is easy to check
$$\ba{ll}
\ns\ds d\Phi^{-1}=\Phi^{-1}\big[\si^2|\Th^*|^2-r-\b\Th^*\big]\Phi^{-1} ds-\Phi^{-1}\si\Th^*dW(s).
\ea
$$
Applying It\^{o}'s formula to $\Phi^{-1}X^*$, we have
\bel{Rewritten-of-X-*}\ba{ll}
\ns\ds  X^*(t)=\Phi(t)x+\Phi(t)\int_0^t\Big[\Phi^{-1}\f^*(\b-\si^2\Th^*)+\Phi^{-1}(l-\si\Th^* h)\Big]ds\\
\ns\ds\qq\qq +\Phi(t)\int_0^t(\Phi^{-1}\si\f^*+\Phi h)dW(s),\ \ t\in[0,T].
\ea
\ee
By the integrability of $(\Th^*,\f^*)$ in (\ref{True-Th-varphi}), for any $p>2$,
\bel{Verification-state-inequality-1}\ba{ll}
\ns\ds  \dbE\Big|\int_0^T|\f^* ||\b -\si^2\Th^* |ds\Big|^{p}\\
\ns\ds \leq \dbE\Big\{\Big[\int_0^T|\f^* |ds\Big]^{\frac p 2} \Big[\int_0^T|\b -\si^2\Th^* |^2ds\Big]^{\frac p 2}\Big\}\\
\ns\ds \leq \Big\{\dbE\Big[\int_0^T|\f^* |ds\Big]^{\frac {p p'} {2}}\Big\}^{\frac {1}{p'}}
 \Big\{\dbE\Big[\int_0^T|\b -\si^2\Th^* |^2ds\Big]^{\frac pq' 2}\Big\}^{\frac {1}{q'}}<\infty,
\ea
\ee
where $p', q'>1$, $\frac{1}{p'}+\frac{1}{q'}=1.$
Similarly, from (H1), for any fixed $p>2$, we get
\bel{Verification-state-inequality-2}\ba{ll}
\ns\ds \dbE\Big|\int_0^T|l-\si\Th^* h|ds\Big|^{p}+\dbE\Big[\int_0^T|\si\f^*+h|^2ds\Big]^{\frac p 2}<\infty.
\ea
\ee
We claim that $\Phi(\cd)=\frac{P_2^*(0)}{P_2^*(\cd)}$ and thus $\Phi$ is bounded. In fact,
$$
\ba{ll}
\ns\ds
\Phi(t)=\exp\Big[\int_0^t(r-\frac{\b\L_2^*}{\si P_2^*}-\frac{|\L_2^*|^2}{2|P_2^*|^2})ds-\int_0^t \frac{\L_2^*}{P_2^*}dW(s)\Big],\  \ t\in[0,T].
\ea
$$
Recall that
$$\ba{ll}
\ns\ds dP_2^*=\Big[-(r-\frac{\b\L_2^*}{\si P_2^*})P_2^*+\frac{|\L_2^*|^2}{P_2^*}\Big]ds+\L_2^*dW(s),
\ea
$$
we get the following via It\^{o}'s formula to $\ln (P_2^*)$,
$$\ba{ll}
\ns\ds d(\ln P_2^*)=(-r+\frac{\b\L_2^*}{\si P_2^*}+\frac{|\L_2^*|^2}{2|P_2^*|^2})ds+\frac{\L_2^*}{P_2^*}dW(s).
\ea
$$
The above expression of $\Phi$ is obvious to see.
Combining this result with (\ref{Rewritten-of-X-*}), (\ref{Verification-state-inequality-1}), (\ref{Verification-state-inequality-2}), we conclude that for any $p>2$, $\dbE\sup\limits_{t\in[0,T]}|X^*(t)|^p<\infty$.

Using again the integrability of $(\Th^*,\f^*)$ in (\ref{True-Th-varphi}), one has
$$\ba{ll}
\ns\ds \dbE\int_0^T|u^*(s)|^2ds\leq \dbE\Big[\sup_{s\in[0,T]}|X^*(s)|^2\int_0^T|\Th^*(s)|^2 ds\Big]+\dbE\int_0^T|\f^*(s)|^2ds<\infty.
\ea
$$

Similarly as above, we have
$$\ba{ll}
\ns\ds \dbE\sup\limits_{t\in[0,T]}|X^{v,\e}_0(t)|^p<\infty,\ \ p>2,\ \ \dbE\int_0^T|u^{v,\e}_0(s)|^2ds<\infty.
\ea
$$
\endpf

\br
Above (H3) is used to verify the pointwise integrability of $\L_2^*$. To explain it via an example, let $r(\cd):=f_1(W(\cd))$, $\frac{\b(\cd)}{\si(\cd)}:=f_2(W(\cd))$, where both function $f_i$ and derivative function $f_i'$ are bounded, $i=1,2.$ In this case, (H3) is obvious to see.
\er

\subsection{Closed-loop equilibrium investment strategies in mean-variance problems}

In this part, we derive the existence of closed-loop equilibrium strategies for mean-variance portfolio selection problems with random coefficients.

To begin with, let us consider
\bel{Equations-for-P-mean-variance-1}\left\{\ba{ll}
\ns\ds d \sP_1^* =- \big[ r  \sP_1^* -\frac{\sL_1^*}{\si \sP_1^*}( \b \sP_1^* +\si\sL_1^*) \big]ds+
\sL_1^* dW(s),\\
\ns\ds d\sP_2^*=-\Big[(r-\frac{\b\sL_1^*}{\si \sP_1^*})\sP_2^*-\frac{\sL_1^*}{\sP_1^*} \sL_2^*\Big]ds+\sL_2^* dW(s),\\
\ns\ds d\sP_3^* =   \frac{(\sP_1^* \b +\sL_1^*\si)(2\sP_1^*\sL_3^*\si+\b\sP_2^*+\si\sL_2^*)}{2\si^2 |\sP_1^*|^2} ds+\sL_3^*dW(s),\\
\ns\ds \sP_1^*(T)=1, \ \ \sP_2^*(T)=-\g, \ \ \sP_3^*(T)=0,
\ea\right.\ee
and define
\bel{Definition-mean-variance-strategy-1}\ba{ll}
\ns\ds  \ \ \Th^*:=-\frac{\sL_1^*}{\si \sP_1^*},\ \  \f^*:=-\frac{2\sP_1^*\sL_3^*\si+\b\sP_2^*+\si\sL^*_2}{2\si^2|\sP_1^*|^2}. %
\ea\ee
The following result is direct conclusion of Theorem \ref{Main-theorem-1}, Theorem \ref{Theorem-main-2}.
\bt
Suppose (H0), (H3) hold. Then there exist $(\sP_i^*,\sL_i^*)$, $i=1,2,3$, satisfying (\ref{Equations-for-P-mean-variance-1}). In addition, above $(\Th^*,\f^*)$ in (\ref{Definition-mean-variance-strategy-1}) is a closed-loop equilibrium operator, and
\bel{Integrability-special-same}\left\{\ba{ll}
\ns\ds \sup_{t\in[\t,T]}\dbE_\t|\Th^*(t)|^2<\infty,\ \ a.s.\ \  \t\in[0,T],\\
\ns\ds  \f^* \in L^p_{\dbF}(\Omega;L^2(0,T;\dbR)), \ \ p>2.
\ea\right.
\ee
\et
\br\label{Remark-mean-variance-1}
If only $r$ is deterministic, then $(\sL_1^*,\sL_2^*)=(0,0)$, and for $s\in[0,T]$,
$$\left\{\ba{ll}
\ns\ds d \sP_i^* =-   r  \sP_i^*  ds,\ \ i:=1,2,\\
\ns\ds d\sP_3^* =  \big[\frac{\b}{\si}\sL_3^*+\frac{\sP_2^*\b^2}{2\si^2\sP_1^*}\big]ds+\sL_3^*dW(s),\\
\ns\ds \sP_1^*(T)=1, \ \ \sP_2^*(T)=-\g, \ \ \sP_3^*(T)=0.
\ea\right.
$$
By the explicit expressions of $\sP_1^*$, $\sP_2^*$, we rewrite $(\sP_3^*,\sL_3^*)$ as follows,
$$\left\{\ba{ll}
\ns\ds d\sP_3^* =  \big[\frac{\b}{\si}\sL_3^*-\frac{\g\b^2}{2\si^2 }\big]ds+\sL_3^*dW(s),\\
\ns\ds  \sP_3^*(T)=0.
\ea\right.
$$
In this case,
$$\Th^*=0,\ \ \f^*=\Big[\frac{\b\g}{2\si^2} -\frac{\sL_3^*}{\si }\Big]e^{-\int_{\cd}^Tr(s)ds}.$$
If $\b$, $\si$ also become deterministic. Then $\sL_3^*=0$, and
$$\f^*=\frac{\g\b }{2|\si|^2}e^{-\int_{\cd}^T r(s)ds}.$$
This result coincides with the analogue in \cite{BC 2010}, \cite{BMZ 2012}, \cite{Huang-Li-Wang-2017}.

We point out several interesting facts which have not been discussed elsewhere.

(1) When $r$ is deterministic, we have $\Th^*=0$, no matter $\b,$ $\si$ are random or not. On the other hand, if
$r$ is random, $\Th^*$ will not degenerate even though $\b$, $\si$ are deterministic. In contrast with $\b$ and $\si$, these conclusions reflect the priority role of $r$ in keeping $\Th^*$.

(2) To make equilibrium strategy rely on initial wealth, the authors introduced state-dependent risk aversion in \cite{BMZ 2012}. Nevertheless, here it is shown that even when the risk aversion is constant, the equilibrium strategy still depends on $X^*$, or initial wealth $x$, as long as $r$ is random.

(3) Suppose $\si$ is random, $\b=0$.  If $r$ is deterministic, then $\f^*=0$, while if $r$ is random, then $\f^*=0$ may not happen. To keep $\f^*$ appear in equilibrium strategy, it indicates that the randomness of $r$ is more important than that of $\si$.

(4) The term $\sL_3^*$ indicates the randomness effects of $\b$, $\si$. This connection is similar as that between $\Th^*$ and $r$ aforementioned.

\er

\br\label{Remark-mean-variance-2}
Let us make a comparison with open-loop equilibrium strategies (\cite{HJZ 2012}, \cite{HJZ 2015}, \cite{Wei-Wang-2017}).

We consider the following systems of equations in $[0,T]$,
$$\left\{\ba{ll}
\ns\ds d\bar\cP_1=- \big[ r \bar\cP_1-(\bar\cP_1\b+\bar\cL_1\si)\frac{\bar\cL_1}{\si\bar \cP_1}\big]ds+
\bar\cL_1dW(s),\\
\ns\ds d\bar\cP_2=-r\bar\cP_2 ds+\bar\cL_2 dW(s),\\
\ns\ds d \bar\cP_3 =-\Big[\big[r-\frac{\b\bar\cL_2}{\si \bar\cP_2}\big]\bar\cP_3-\frac{\bar\cL_2}{\bar \cP_2}\bar\cL_3\Big] ds+\bar\cL_3dW(s),\\
\ns\ds d\bar\cP_4=- \big[\bar\cP_1\b+\bar\cL_1\si\big]\Big\{\frac{\b \bar\cP_3+\si \bar\cL_3}{\si^2\bar\cP_2\bar\cP_1}-\frac{\bar\cL_4}{\si \bar \cP_1}\Big\} ds+\bar\cL_4dW(s),\\
\ns\ds \bar\cP_1(T)=1,\ \ \bar\cP_2(T)=-2, \ \ \bar\cP_3(T)=-\g, \ \ \bar\cP_4(T)=0,\ \
\ea\right.$$
and define
$$\ba{ll}
\ns\ds  \bar\Th:=-\frac{\bar\cL_1}{\si \bar\cP_1},\ \ \  \bar\f:=\frac{\b \bar\cP_3+\si \bar\cL_3}{\si^2 \bar\cP_2\bar\cP_1}-\frac{\bar\cL_4}{\si \bar\cP_1}.\ \
\ea
$$
For $\nu,s\in[0,T]$, suppose $D_\nu r(s)$ exists and $\big|D_\nu r(s)\big|\leq K$. Using Proposition 3.9 in \cite{Wei-Wang-2017}, we see that $\bar u:=\bar\Th \bar X+\bar\f$ is an open-loop equilibrium strategy, and $(\bar\Th,\bar\f)$ is a pair of open-loop equilibrium operator which is independent of initial wealth.

If $r$ is deterministic, then $\bar\cL_i=0,$ $i=1,2,3$,
$$\left\{\ba{ll}
\ns\ds d\bar\cP_i=-   r \bar\cP_i ds, \ \ i:=1,2,3,\\
\ns\ds d\bar\cP_4=-  \bar\cP_1\b \Big\{\frac{\b \bar\cP_3 }{\si^2\bar\cP_2\bar\cP_1}-\frac{\bar\cL_4}{\si \bar \cP_1}\Big\} ds+\bar\cL_4dW(s),\\
\ns\ds \bar\cP_1(T)=1,\ \ \bar\cP_2(T)=-2, \ \ \bar\cP_3(T)=-\g, \ \ \bar\cP_4(T)=0.
\ea\right.$$
By the explicit expression of $\bar\cP_i$, $i=1,2,3$, we can rewrite $(\bar\cP_4,\bar\cL_4)$ as,
$$\left\{\ba{ll}
\ns\ds d\bar\cP_4=\big[\frac{\b}{\si }\bar\cL_4-\frac{\b^2\g}{2\si^2}\big]ds+\bar\cL_4dW(s),\\
\ns\ds \bar\cP_4(T)=0.
\ea\right.$$
In this case,
$$\bar\Th=0,\ \ \bar\f=\Big[\frac{\b\g}{2\si^2}-\frac{\bar\cL_4}{\si\bar\cP_1}\Big]e^{-\int_{\cd}^Tr(s)ds}.$$

We make a few interesting points, which have not mentioned in existing papers to our best.

(1) As to $\Th^*$, $\bar\Th$, they are equal if they exist, even when $(r,\b,\si)$ are random. In other words, both closed-loop and open-loop equilibrium strategies have the same manner in depending on equilibrium wealth process.

(2) If $r$ is deterministic, we have further equality between $\f^*$, $\bar\f$, even when $(\b,\si)$ are random. This means that mean-variance portfolio selection problem admits a pair of closed-loop, open-loop equilibrium strategies that are the same with each other.

(3) Above two equality conclusions indicate the peculiar role with constant risk averison, since they do not happen even in the same framework with state-dependent risk aversion (\cite{HJZ 2012}).

\er


\subsection{Closed-loop equilibrium hedging strategies in variance hedging problems}

In this part, we give the explicit closed-loop equilibrium strategies of dynamic variance hedging problems in non-Markovian setting, which has not been done in the literature to our best.

Since $\dbE|\xi|^{k}<\infty$, $k>2$, we have
$$\ba{ll}
\ns\ds \dbE\sup_{s\in[0,T]}|\l(s)|^{k}<\infty,\ \  \dbE\Big[\int_0^T|\zeta(s)|^2ds\Big]^{\frac k 2}<\infty.
\ea
$$

We consider the following system of equations
\bel{Equations-for-P-variance-hedging-1}\left\{\ba{ll}
\ns\ds d \sP_1^* =- \big[ r  \sP_1^* -\frac{\sL_1^*\b}{\si }-\frac{|\sL_1^*|^2}{\sP_1^*}  \big]ds+
\sL_1^* dW(s),\\
\ns\ds d\sP_2^* =  \Big[\sL_2^*(\frac{\b}{\si}+\frac{\sL_1^*}{\sP_1^*}) -\sP_1^* (r \l+\frac{\b}{\si}\zeta)\Big] ds+\sL_2^*dW(s),\\
\ns\ds \sP_1^*(T)=1, \ \ \sP_2^*(T)=0,
\ea\right.\ee
and define
\bel{Definition-variance-hedging-strategy-1}\ba{ll}
\ns\ds  \ \ \Th^*:=-\frac{\sL_1^*}{\si \sP_1^*},\ \  \f^*:=
-\frac{ \sL_2^*}{\si \sP_1^*}+\frac{\zeta}{\si }.  %
\ea\ee

The following result is direct conclusion of Theorem \ref{Main-theorem-1}, Theorem \ref{Theorem-main-2}.
\bt
Suppose (H0), (H3) hold, $\dbE|\xi|^k<\infty$, $k>2$. Then there exist $(\sP_i^*,\sL_i^*)$, $i=1,2$, satisfying (\ref{Equations-for-P-variance-hedging-1}). In addition, above $(\Th^*,\f^*)$ in (\ref{Definition-variance-hedging-strategy-1}) is a closed-loop equilibrium operator, $\pi^*:=\Th^* X^*+\f^*$ is a closed-loop equilibrium  hedging strategy, and (\ref{Integrability-special-same}) holds with $p$ replaced by $k$.
\et
\br
As one part of equilibrium operator, above $\Th^*$ is the same as the analogue of equilibrium operator of mean-variance portfolio selection problems. Therefore, if $r$ is deterministic, and $\b,\si$ are random, one still has $\Th^*=0$, which means the hedging strategy $\pi^*$ becomes wealth independent.
\er
\br
If $r$ is deterministic, $\b=0$. Then $\sL^*_1=0$, and for $s\in[0,T]$,
$$\left\{\ba{ll}
\ns\ds d\sP_2^* = -\sP_1^* r \l ds+\sL_2^*dW(s),\\
\ns\ds \sP_2^*(T)=0.
\ea\right.
$$
The randomness of $\l$ implies that $\f^*$ may not degenerate. This is different from the corresponding study in mean-variance portfolio selection problems (see Remark \ref{Remark-mean-variance-1}).
\er

\br
The dynamic variance hedging problem with different formulation was discussed in \cite{Basak-Chabakauri-2012} with deterministic coefficients. Our study can be adjusted into their situation with general random coefficients. We hope to demonstrate more related results in future.
\er

\section{Concluding remark}

In this paper, we discuss the dynamic mean-variance portfolio selection problem and dynamic variance hedging problem in a unified manner. The involved coefficients are allowed to be random, which implies the introducing of systems of backward stochastic Riccati equations. Due to the non-Markovian setting, several interesting phenomena, which are concealed in existing literature with deterministic coefficients, are revealed here for the first time.

Without further essential difficulties, our investigation also works when investment strategy, hedging strategy, as well as Brownian motion are multi-dimensional. The uniqueness of closed-loop equilibrium strategies of these financial problems are still under consideration. We hope to present it in future publications.

\end{document}